\documentclass[12pt]{article}
\usepackage{amssymb,latexsym,amsmath,graphics}
 
\newtheorem{defn}{Definition}[section]
\newtheorem{lemma}[defn]{Lemma}

\newtheorem{corollary}[defn]{Corollary}

\DeclareMathOperator{\eval}{eval}

\def\endproof{\relax\ifmmode\expandafter\endproofmath\else
 \unskip\nobreak\hfil\penalty50\hskip.75em\hbox{}\nobreak\hfil\bull
  {\parfillskip=0pt \finalhyphendemerits=0 \bigbreak}\fi}

\def\bull{\vbox{\hrule\hbox{\vrule
    \kern3pt\vbox{\kern6pt}\kern3pt\vrule}\hrule}}

\makeatletter
\newfont{\footsc}{cmcsc10 at 8truept}
\newfont{\footbf}{cmbx10 at 8truept}
\newfont{\footrm}{cmr10 at 10truept}
\renewcommand{\ps@plain}{%
\renewcommand{\@oddfoot}{\footsc the electronic journal of combinatorics
  {\footbf 12} (2005), \#R00\hfil\footrm\thepage}}
\makeatother 

\def\dsp{\displaystyle}
\def\S{\sum \limits_{n= 0}^\infty}

\makeatletter 
  \def\@seccntformat#1{\csname the#1\endcsname.\quad} 
\makeatother 

\let\oldsection=\section
\def\section#1{\oldsection{\kern -13pt #1}}

\let\oldsubsection=\subsection
\def\subsection#1{\oldsubsection{\kern -11pt #1}}

\title{A  triple lacunary generating function for Hermite polynomials}

\author{Ira M. Gessel\thanks
   {Partially supported by NSF Grant  DMS-0200596}\\
\small Department of Mathematics\\[-0.8ex]
\small Brandeis University, Waltham, MA, USA\\[-0.8ex]
\small \texttt{gessel@brandeis.edu}
\and
Pallavi Jayawant\\
\small Department of Mathematics\\[-0.8ex]
\small Bates College, Lewiston, ME, USA\\[-0.8ex]
\small \texttt{pjayawan@bates.edu}}

\date{\small Submitted: January 21, 2004;  Accepted: March 29, 2005 ;
 Published: \\
\small Mathematics Subject Classifications: 05A15, 05A19, 05A40, 33C45}

\begin{document}

\maketitle

\begin{abstract}
Some of the classical orthogonal polynomials such as Hermite, Laguerre,
Charlier, etc. have been shown to be the generating polynomials for certain
combinatorial objects. These combinatorial interpretations are used to prove
new identities and generating functions involving these polynomials.  In
this paper we apply Foata's approach to generating functions for the Hermite
polynomials to obtain a triple lacunary generating function. We define
renormalized Hermite polynomials $h_n(u)$ by
$$\S h_n(u) \frac{z^n}{ n!}=e^{uz+{z^2\!/2}}.$$
and give a combinatorial proof of the following  generating
function:
\begin{equation*}
\sum \limits_{n= 0}^\infty h_{3n}(u) \dsp{\frac{z^n}{n!}}= \frac
{e^{(w-u)(3u-w)/6}}{\sqrt{1-6wz}}
\sum\limits_{n= 0}^\infty\dsp{\frac
{(6n)!}{2^{3n}(3n)!(1-6wz)^{3n}}
\frac{z^{2n}}{(2n)!}},
\end{equation*}
 where
$w=(1-\sqrt{1-12uz})/6z=uC(3uz)$ and $C(x)=(1-\sqrt{1-4x})/(2x)$
is the Catalan generating function.
\noindent
We also give an
umbral proof of this generating function.
\end{abstract}

\section{Introduction} \label{S:intro}

The Hermite polynomials $H_n(u)$ may be
defined by the exponential generating function
\begin{equation}
\S H_n(u) \frac {z^n} {n!}=e^{2uz-z^2}. 
\end{equation}
In this paper we will give a combinatorial proof of an
identity for Hermite polynomials. For our combinatorial
interpretation, it is more convenient to take  a different
normalization of the Hermite polynomials, which makes all the
coefficients positive. Therefore, we work with the
polynomials $h_n(u)=\dsp{\frac{i^n}{2^{n/2}} H_n
\biggl(\frac{-iu}{\sqrt 2}\biggr)}$, where
$i=\sqrt{-1}$, which have the generating function
$$\S h_n(u) \frac{z^n}{n!}=e^{uz+{z^2\!/2}}.$$
All of our formulas for $h_n(u)$ are easily converted into
formulas for $H_n(u)$.

Foata \cite{foata:doetsch} gave a combinatorial proof of  Doetsch's identity
\cite{doetsch}
giving a generating function
for $h_{2n}(u)$:
\begin{equation} \label{E:doetsch:evenmodify}
\S h_{2n}(u) \frac {z^n} {n!} =(1-2z)^{-1/2}  \exp{\biggl(\frac
{u^2z}{1-2z}\biggr)}.
\end{equation}
 We will prove the following generating function for $h_{3n}(u)$:

\begin{equation}\label{E:mainf}\sum \limits_{n= 0}^\infty h_{3n}(u) \dsp\frac{z^n}{n!}= \frac
{e^{(w-u)(3u-w)/6}}{\sqrt{1-6wz}}
\sum\limits_{n= 0}^\infty\dsp{\frac
{(6n)!}{2^{3n}(3n)!(1-6wz)^{3n}}
\frac{z^{2n}}{(2n)!}},\end{equation}
in which
$w=(1-\sqrt{1-12uz})/6z=uC(3uz)$, where  $C(x)=(1-\sqrt{1-4x})/2x$
is the Catalan number generating function.  Note that the formula can be
written in terms of hypergeometric series:
$$\sum \limits_{n= 0}^\infty h_{3n}(u) \dsp\frac{z^n}{n!}= {\frac
{e^{(w-u)(3u-w)/6}}{\sqrt{1-6wz}}} \,
_2F_0\biggl(\frac 1 6,\frac 5
6;-;\frac{54z^2}{(1-6wz)^3}\biggr),$$ where $\dsp{_2F_0(a, b, -;
z) = \sum \limits_{n \ge 0} (a)_n \, (b)_n \frac{z^n}  {n \, !}}$,
and
$${(\alpha)_n = \alpha(\alpha+1)\ldots (\alpha+n-1)}$$ is the rising factorial.

We prove formula (\ref{E:mainf}) by two methods---umbral and combinatorial. In
section \ref{S:umbralproof}, we define an umbra and study some of its properties and then give the umbral proof. An umbral proof of a generating function for
$h_{2m+n}(u)$ is given in \cite{gessel}.

In section \ref{S:combpf} we prove (\ref{E:mainf}) combinatorially by showing that both sides enumerate the
same weighted objects. We first describe the combinatorial interpretation of the Hermite polynomials and then give the details of the weighted objects counted by both sides of the formula. 
 
By using these methods, it would be possible to give umbral and combinatorial proofs of a more general generating function for $h_{3m+2n+k}(u)$.

\section{The Umbral Proof}\label{S:umbralproof}

Rota and Taylor in \cite{rota-taylor:intro} laid a rigorous foundation for
the classical umbral calculus. They consider a vector space of polynomials in
several variables or ``umbrae" and define the linear functional $\eval$ on
it. A sequence $(a_n)$ is represented by an umbra $A$ if $\eval(A^n)=a_n$
for all $n$. In practice, the word $\eval$ is usually dropped and we simply
write $A^n= a_n$ with the understanding that the functional has been applied.
When we write $f(A)=g(A)$, we mean $\eval(f(A))=\eval(g(A))$. We consider
formal power series $f(t)$ with coefficients in a ring of formal power series
$R={\mathbf {Q}}[[x, y, z, \cdots]]$.
\begin{defn} A formal power series $f(t)=\S f_nt^n$ is {\it admissible} if for
every monomial $x^iy^jz^k \cdots $  in $R$, the coefficient of
$x^iy^jz^k \cdots $ in $f_n$ is nonzero for only finitely many values of $n$.
\end{defn} So, for example, $f(t)=e^{xt}$ is an admissible formal power series,
while $f(t)=e^t$ is not. Some computations similar to those of this section can
be found in section 4 of Gessel \cite{gessel}.

\subsection{The umbra $M$ and its properties} \label{S:umbra}

We now define an umbra $M$ whose relation to the Hermite polynomials will be
described in the next section. In this section, we study some properties of
this umbra that we will need. We define the umbra $M$ by $e^{Mz} =
e^{z^2\!/2}$ so that
\begin{equation} \label{E:mdefn}
M^n= \begin{cases}
\dsp\frac{(2k)!}{2^k k!}, &\text{if $n=2k$}\\
0, &\text{if $n$ is odd.}
\end{cases}
\end{equation}

The following two formulas hold for $M$.
\begin{lemma} \label{L:fm}

(i) For any admissible formal power series $f$, $e^{Mz}
f(M)=e^{z^2\!/2}f(M+z).$

(ii) $e^{M^2z}=\dsp{\frac {1}{\sqrt{1-2z}}}.$
\end{lemma}

\begin{proof}
(i) We first observe that it is sufficient to prove that the formula holds
for $f(t)=e^{tx}$: If the formula is true for $f(t)=e^{tx}$, then comparing
coefficients of $x^n/n!$ on both sides we get that the formula holds for
$f(t)=t^n$. But this implies by linearity that the formula is true for all
admissible formal power series $f$. To prove the formula for $f(t)=e^{tx}$,
we have

$$e^{Mz}e^{Mx}=e^{M(z+x)}=e^{z^2\!/2+zx+x^2/2}=e^{z^2\!/2}e^{zx}e^{Mx}=e^{z^2\!/2}e^{(M+z)x}.$$
(ii) We have
$$e^{M^2z}=\sum\limits_{n= 0}^\infty M^{2n}\frac{z^n}{n!}=\sum\limits_{n= 0}^\infty
\frac{(2n)!}{2^n n!}\frac{z^n}{n!}=\frac {1}{\sqrt{1-2z}}.$$
\end{proof}
\begin{corollary}\label{C:msquare}
For any admissible formal power series $f$,
$$e^{M^2z}
f(M)=\dsp{\frac {1}{\sqrt{1-2z}}f\biggl(\frac
{M}{\sqrt{1-2z}}\biggr)}.$$
\end{corollary}
\begin{proof}
As in the lemma, it is sufficient to prove the formula for
$f(t)=e^{tx}$. Here $e^{M^2z} f(M)=e^{M^2z+Mx}.$ If we apply Lemma
\ref{L:fm} (i)
directly, we cannot
eliminate the linear term in $M$.
So
we introduce a parameter $\alpha$ and rewrite $e^{M^2z+Mx}$ as
$e^{M\alpha}e^{M^2z+M(x-\alpha)}$. We will choose a value for $\alpha$ later. Now
applying Lemma \ref{L:fm} (i) we get
$$e^{M\alpha}e^{M^2z+M(x-\alpha)}=e^{\alpha^2/2}e^{(M+\alpha)^2z+(M+\alpha)(x-\alpha)}
=e^{M^2z+M(x-(1-2z)\alpha)+z\alpha^2+x\alpha-\alpha^2/2}.$$
Now we choose the value of $\alpha$ to eliminate the term in $M$ on the
right; i.e., we solve $x-(1-2z)\alpha=0$ and get $\alpha=x/(1-2z)$. Substituting
this value of $\alpha$ in the above expression and simplifying, we
obtain that $e^{M^2z+Mx}=e^{M^2z} e^{x^2/2(1-2z)}$. By Lemma
\ref{L:fm} (ii),
this is equal to $\dsp{\frac {1}{\sqrt{1-2z}}}
e^{x^2/2(1-2z)}$. But applying the definition of $M$ directly gives
$$\dsp{\frac {1}{\sqrt{1-2z}}f\biggl(\frac{M}{\sqrt{1-2z}}\biggr)}
  = \frac {1}{\sqrt{1-2z}} e^{Mx/\sqrt{1-2z}}=\frac {1}{\sqrt{1-2z}} e^{x^2/2(1-2z)}.
$$
\end{proof}
Applying the corollary to $f(M)=e^{M^3x}$ gives the formula
\begin{equation}\label{E:cor}e^{M^2z+M^3x}=\frac{\exp\dsp{\biggl({\frac {M^3x}
{(1-2z)^{3/2}}}\biggr)}}{\sqrt{1-2z}},
\end{equation}
which we will need in the next section.

\subsection{Proof of the formula} \label{S:umbralpf}

To prove formula (\ref{E:mainf}), we first express the
 Hermite polynomial $h_n(u)$ in terms of the
umbra $M$. We have
$$\S h_n(u) \frac{z^n}{n!}=e^{uz+{z^2\!/2}}= e^{(u+M)z}.$$
Comparing the coefficients of $z^n/n!$ on both sides, we get
$h_n(u)=(u+M)^n$. Using this, we get
$$\S h_{3n}(u) \frac{z^n}{n!}=\S (u+M)^{3n} \frac{z^n}{n!}
=e^{(u+M)^3z}=e^{(u^3+3Mu^2+3M^2u+M^3)z}.$$ We follow the same
procedure as in the proof of Corollary \ref{C:msquare}. We
introduce a parameter $\alpha$ and rewrite the last expression as
$e^{M\alpha}e^{u^3z+M(3u^2z-\alpha)+3M^2uz+M^3z}$. Now applying Lemma
\ref{L:fm} (i) we get
\begin{multline}
e^{M\alpha}e^{u^3z+M(3u^2z-\alpha)+3M^2uz+M^3z}=e^{\alpha^2/2}e^{u^3z+(M+\alpha)(3u^2z-\alpha)
    +3(M+\alpha)^2uz+(M+\alpha)^3z}  \\
=e^{z\alpha^3+3uz\alpha^2-\alpha^2/2+3u^2z\alpha+u^3z+M(3z\alpha^2
    +(6uz-1)\alpha+3u^2z)+M^2(3z\alpha+3uz)+M^3z}\label{mexp}.\end{multline}
In order to apply formula (\ref{E:cor}), we need to eliminate the
linear term in $M$. So we choose a value of $\alpha$ that makes the
coefficient of $M$ equal to zero. By solving a quadratic and taking
the solution with a power series expansion,  we get
$\alpha=(1-\sqrt{1-12uz})/6z-u.$  We know that the
Catalan generating function $C(x)$ is given by
$(1-\sqrt{1-4x})/2x$. In terms of this generating function,
$\alpha=uC(3uz)-u=w-u$ where $w=uC(3uz)$. Using this value of $\alpha$ to
simplify expression (\ref{mexp}), we get that
$$\S h_{3n}(u) \frac{z^n}{n!}=\exp(w^3z-(w-u)^2/2+3M^2wz+M^3z).$$

We know that $C(x)$ satisfies the equation $C(x)=1+x(C(x))^2$.
Substituting $x=3uz$ and using the fact that $C(3uz)=w/u$, we
obtain $w=u+3w^2z$; i.e., $w-u=3w^2z$. This gives us
$$\S h_{3n}(u) \frac{z^n}{n!}=e^{(w-u)(3u-w)\,/6}e^{3M^2wz+M^3z}.$$
Applying formula (\ref{E:cor}) with $f(t)=t^3z$, we get
that
$$\sum \limits_{n= 0}^\infty h_{3n}(u) \dsp\frac{z^n}{n!}= \frac
{e^{(w-u)(3u-w)/6}}{\sqrt{1-6wz}}
\exp{\biggl(\frac{M^3z}{(1-6wz)^{3/2}}\biggr)}.$$ Then writing the
second exponential function as a series and using (\ref{E:mdefn}),
we obtain the final result:
$$\sum \limits_{n= 0}^\infty h_{3n}(u) \dsp\frac{z^n}{n!}= \frac
{e^{(w-u)(3u-w)/6}}{\sqrt{1-6wz}}
\sum\limits_{n= 0}^\infty\dsp{\frac
{(6n)!}{2^{3n}(3n)!(1-6wz)^{3n}} \frac{z^{2n}}{(2n)!}}.$$

 Now we turn to the combinatorial method of proof. We begin with a
 combinatorial interpretation of the Hermite polynomials that will
 be used in the combinatorial proofs.

 
 \section{The Combinatorial Proof}
 \label{S:combpf}

We assume that the reader is familiar with enumerative
applications of exponential generating functions, as described,
for example, in \cite[Chapter 5]{stanley:vol2} and \cite{bergeron}. 
The product formula and the exponential formula for exponential
generating functions discussed in these references play an
important role in the combinatorial proofs.

\subsection{Combinatorial interpretation of Hermite polynomials} 
\label{S:combint}

The exponential generating function $uz+z^2/2$ counts
sets with one or two elements, where a one-element set
is weighted $u$. Then by the exponential formula, the
coefficient of $z^n/n!$ in $e^{uz+z^2/2}$, which is 
$h_n(u)$, is the generating polynomial for partitions
of an $n$-element set into blocks of size one  or two,
where each block of size one is weighted $u$. (If we 
used $H_n(u)$ as Foata did,  instead of $h_n(u)$, we
would need to attach a weight of $-2$ to each
two-element block and a weight of $2u$ to each one
element block.) It is convenient to represent these
partitions as graphs in which the vertices in a
two-element block are joined by an edge. We call these
graphs, in which every vertex has degree at most one,
{\it matchings\/}.

Thus the Hermite polynomial $h_n(u)$ can be viewed as the
generating polynomial for the number of vertices of degree zero
over the set of all matchings on $n$ vertices, where each vertex
of degree zero is assigned the weight $u$. With this combinatorial
interpretation we will prove the following formula:

\begin{equation}\label{E:mainform}\sum \limits_{n= 0}^\infty h_{3n}(u) \dsp\frac{z^n}{n!}=
e^{(w-u)(3u-w)/6} \frac {1} {\sqrt{1-6wz}}
\sum\limits_{n= 0}^\infty\dsp{\frac
{(6n)!}{2^{3n}(3n)!\,(1-6wz)^{3n}}
\frac{z^{2n}}{(2n)!}}.\end{equation}

We will describe the graphs enumerated by the left side of
formula (\ref{E:mainform}). Then we will describe the same graphs in terms of their
connected components and use the product formula for exponential
generating functions to complete the proof.

\subsection{Graphs counted by the left side.}
\label{S:leftside}

\begin{figure}[htbp]
\begin{center}
\scalebox{.70}{\includegraphics{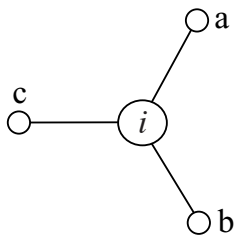}} 
\caption{A label and three marked vertices} \label{Fi:trivunit}
\end{center}
\end{figure}

\begin{figure}[htbp]
\begin{center}
\scalebox{.63}{\includegraphics{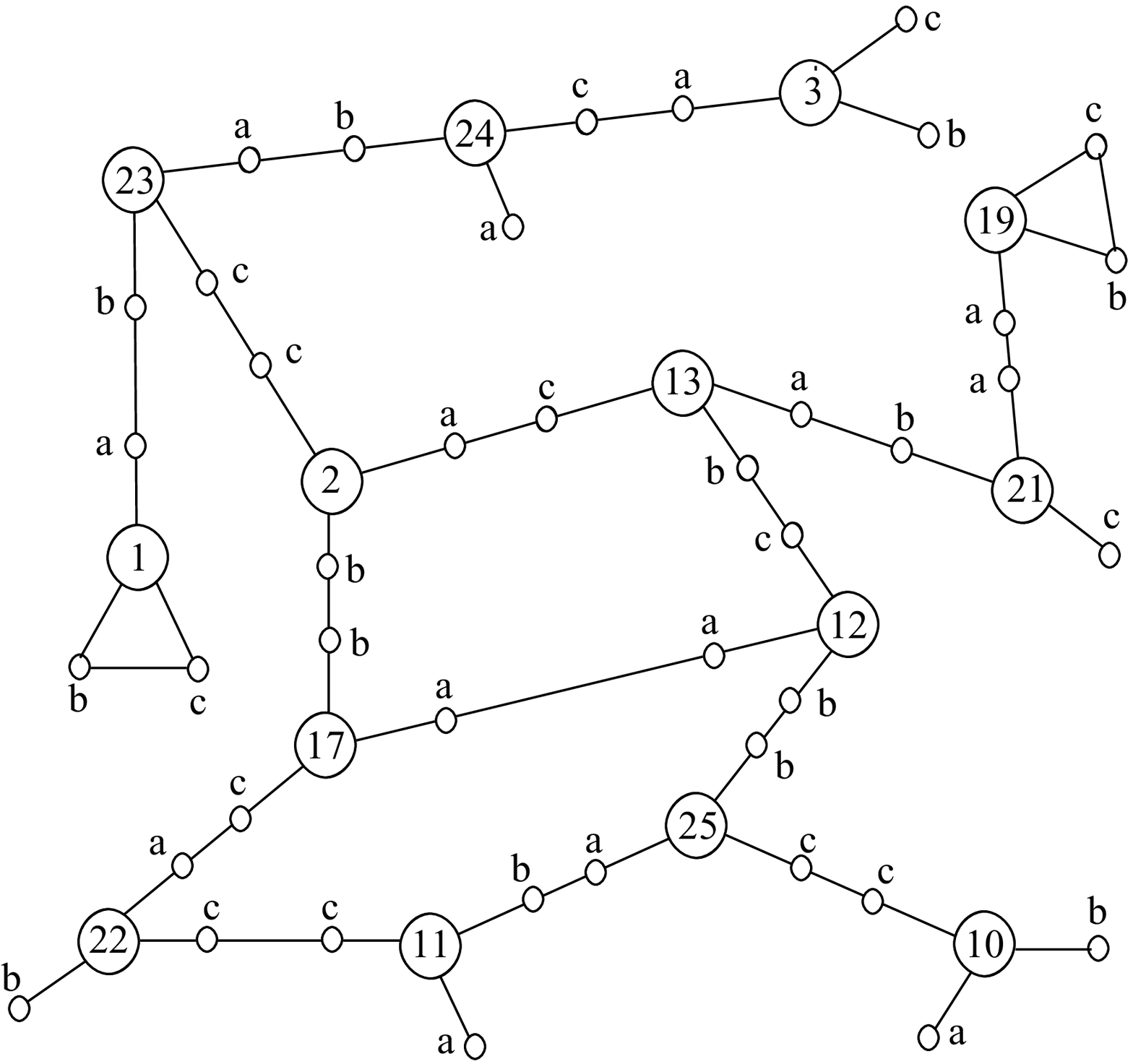}}
\scalebox{.63}{\includegraphics{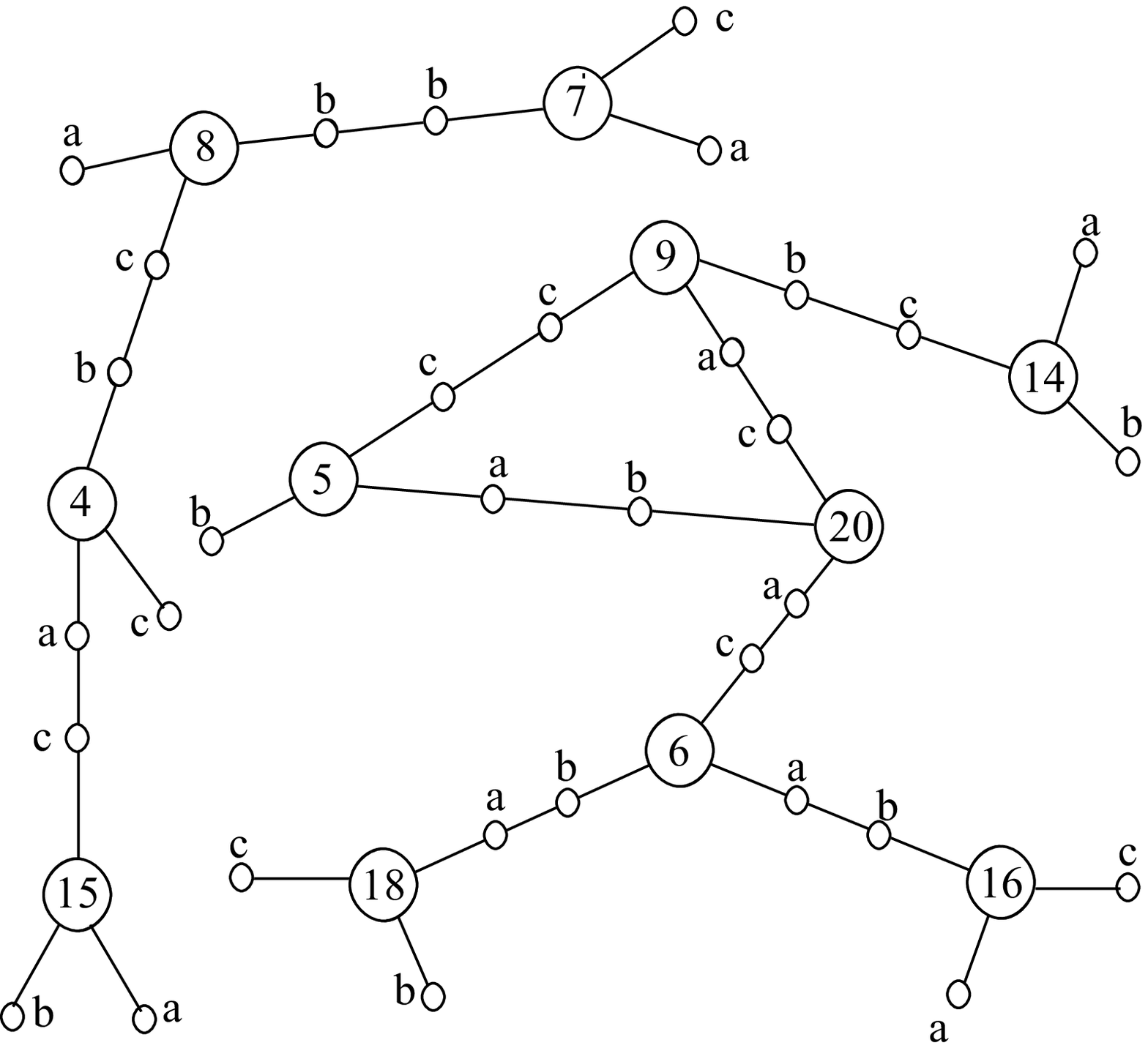}}
 \caption{Graph counted by $h_{75}(u) z^{25}/25!$} \label{Fi:originalgraph}
\end{center}
\end{figure}

In order to give a combinatorial interpretation to  $\sum_{n=0}^\infty
h_{3n}(u)x^n/n!$, we must interpret $h_{3n}(u)$, which counts matchings of a
$3n$-element set, as counting labeled  objects with $n$ labels. To
accomplish this, we take $n$ labels and attach to each one  three vertices
marked $a$, $b$, and $c$. 
Figure \ref{Fi:trivunit} shows a labeled vertex connected to three marked
vertices.
Then $h_{3n}(u)$ counts graphs constructed by taking $n$ components like that in
Figure~\ref{Fi:trivunit}, with labels 1, 2, \dots, $n$ and adding to them a matching of
the $3n$ marked vertices, where each unmatched marked vertex is weighted $u$.
Figure~\ref{Fi:originalgraph} shows such a graph.

Let $G$ be the set of all graphs enumerated by the left side of formula
(\ref{E:mainform}).  For a graph in $G$, each vertex has degree one, two, or
three.  The trivalent vertices are labeled with the integers 1 to $n$, where
$n$ is the  total number of trivalent vertices, and are weighted by the
exponential generating function variable $z$. The bivalent and monovalent
vertices are marked with $a$, $b$, and $c$; the monovalent vertices have
weight $u$ and the bivalent vertices have weight 1.

We first make a preliminary simplification of the graphs in $G$: we eliminate
all the bivalent marked vertices, moving their ``marks" to the adjacent
trivalent vertices. 

\begin{figure}[htbp]
\begin{center}
\scalebox{.70}{\includegraphics{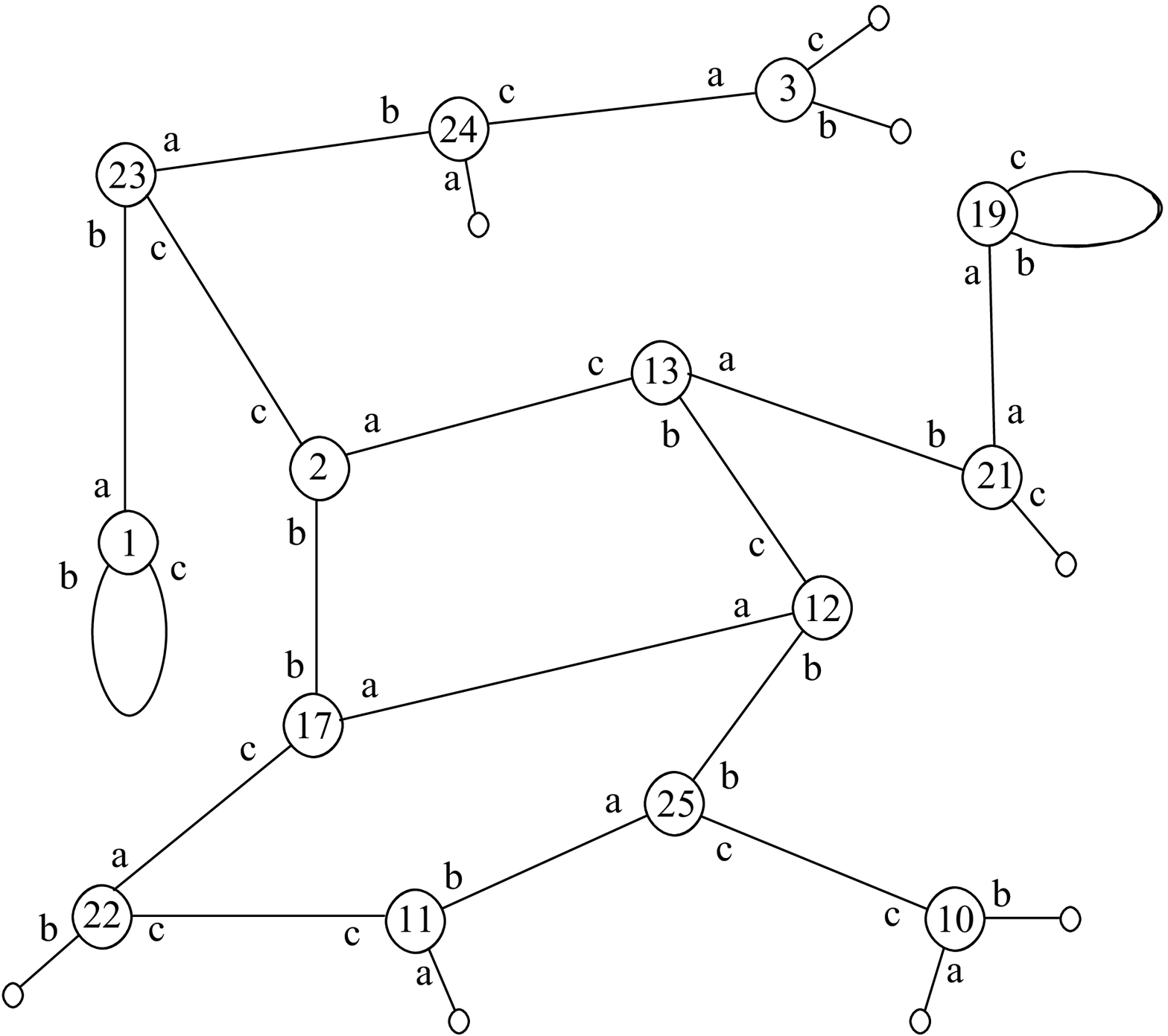}}
\scalebox{.70}{\includegraphics{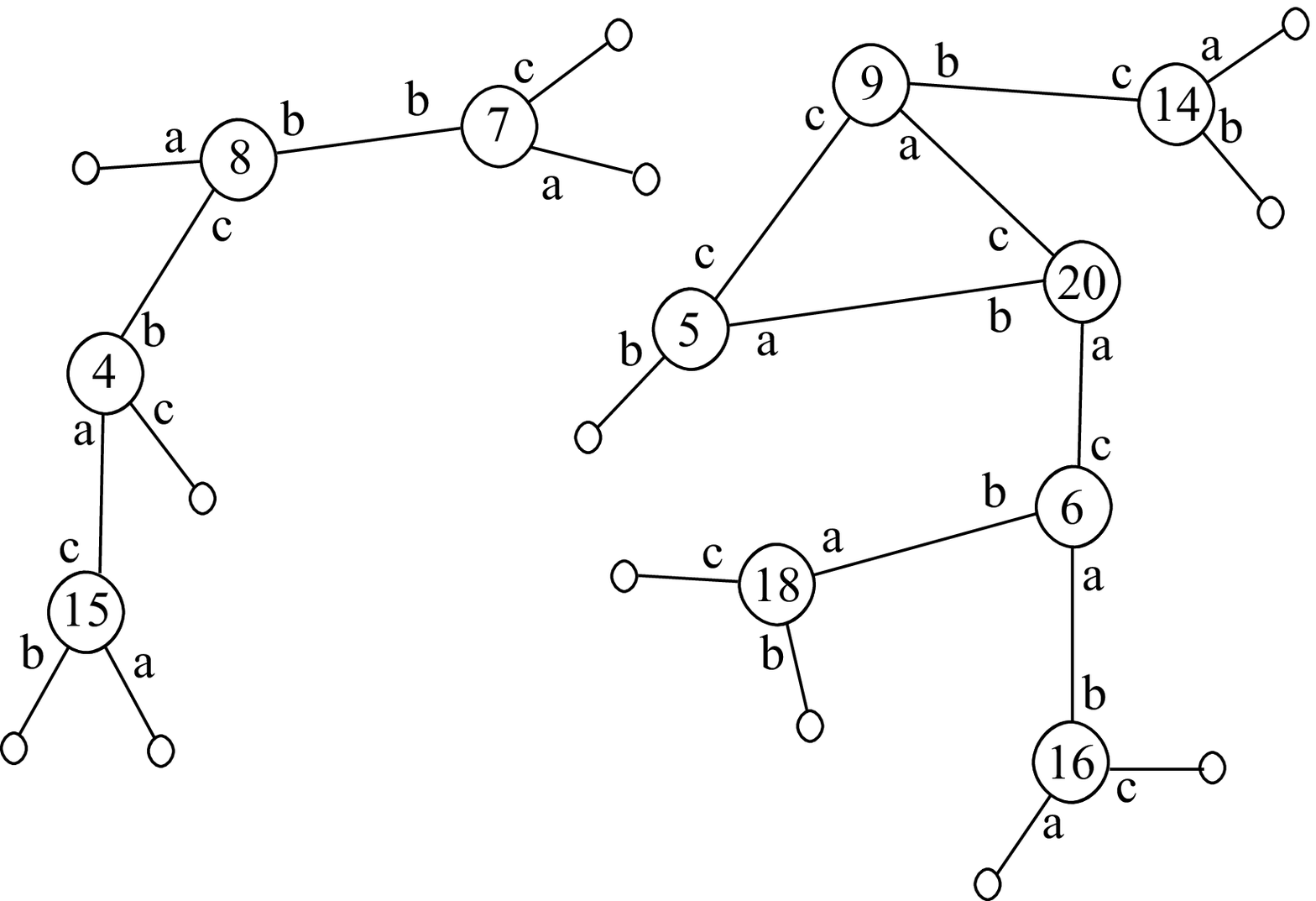}}
\caption{Graph in Figure \ref{Fi:originalgraph} with bivalent
vertices eliminated} \label{Fi:bivertelim}
\end{center}
\end{figure}

More precisely, we think of each edge joining two trivalent vertices as
consisting of two ``half-edges", each of which has a mark.  Although we retain
the monovalent vertices, we move their marks to the  half-edge  of the adjacent
trivalent vertices. Note that some trivalent vertices now have loops. Figure
\ref{Fi:bivertelim}, which shows the transformed version of the graph in Figure
\ref{Fi:originalgraph}, should make this simplification clear. 

We will express the generating function for graphs in $G$ as a product of 
three factors corresponding to connected components  with no cycles, with
exactly one cycle,  and with at least two cycles. Note that a loop is a cycle.

\subsection{Counting $w$-trees}
The graph in
Figure \ref{Fi:bivertelim} has three connected components. One component is a
tree. The component with one cycle may be viewed as a cycle of trees. The third
component may also be decomposed into trees and cycles.  Thus the first step of
our enumeration is to count the rooted trees from which we will construct our
graphs, which we call \emph{$w$-trees}.

\begin{defn}
A $w$-tree is a rooted tree with labeled trivalent internal vertices and
unlabeled leaves in which the three half-edges incident with every internal
vertex are  marked with the letters $a$, $b$ and
$c$, and one of the half-edges incident with the root has no matching half-edge.
\end{defn}

\begin{figure}[htbp]
\begin{center}
\scalebox{.70}{\includegraphics{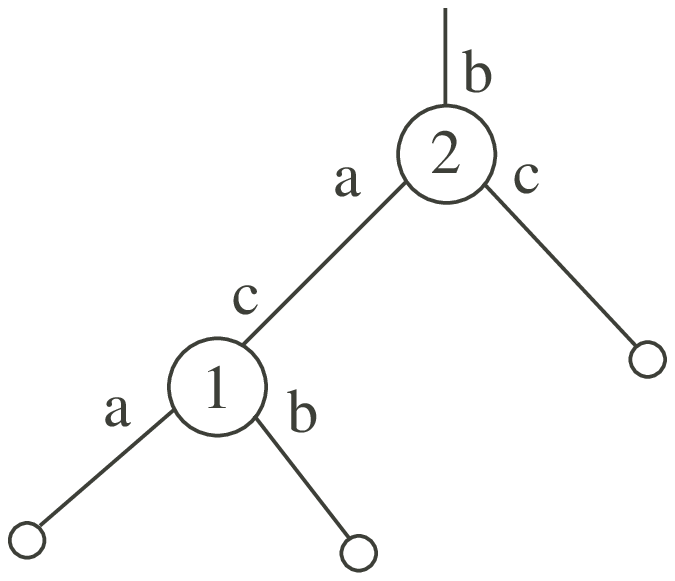}} \caption{$w$-tree}
\label{Fi:wt}
\end{center}
\end{figure}

Figure \ref{Fi:wt} shows a $w$-tree.
A special case of a \mbox{$w$-tree} is the one with no internal
vertices and no marks: it is simply a leaf with a half-edge attached to it.

We denote by $w$ the exponential generating function in $z$ for $w$-trees where
each internal vertex is weighted with $z$ and each leaf is weighted with $u$. 
Let $W_n$ be the number of $w$-trees with $n$ internal vertices. There are
$2^n$ ways to draw such a tree (with the root at the top) by choosing an
ordering of the two children of each internal vertex, so there are $2^nW_n$
such drawings. But these drawings can be counted in another way. If we remove
the labels and marks from such a drawing, we have a binary tree counted by the
Catalan number $C_n = \frac1{n+1}\binom{2n}n$. The marks can be added in $6^n$
ways and the labels in $n!$ ways. Thus $2^nW_n=6^n n!\, C_n$, so $W_n = 3^n
n!\, C_n$. Each such tree has $n+1$ leaves, so the exponential generating
function for weighted $w$-trees is 
$w=\sum_{n=0}^\infty 3^n n!\, C_n u^{n+1} z^n/n! = uC(3uz)$. This relation to the Catalan generating function gives the equation 
$w=u+3w^2z$ as shown in section \ref{S:umbralpf}.

\subsection{Acyclic components.}
The connected components of graphs in $G$ with no cycles are trees.  Such a
tree appears in the lower left corner of Figure \ref{Fi:bivertelim}. 

To count these trees we use a known method (see, e.g., \cite{flajolet}) for relating labeled unrooted trees of various types to rooted trees: we subtract the edge-rooted versions from the vertex-rooted versions.

First we count trees rooted at a trivalent vertex. In such a tree,  the three half-edges at the root are joined to $w$-trees.
So we can construct such a rooted tree by taking an ordered triple of
$w$-trees, with exponential generating function $w^3$,  and a new root vertex with its three half-edges marked $a$, $b$, and $c$, and attaching the first $w$-tree to half-edge $a$, the second to half-edge $b$, and the third to half-edge $c$. Thus the exponential generating function for these vertex-rooted trees is
$w^3z$.

\begin{figure}[htbp]
\begin{center}
\scalebox{.80}{\includegraphics{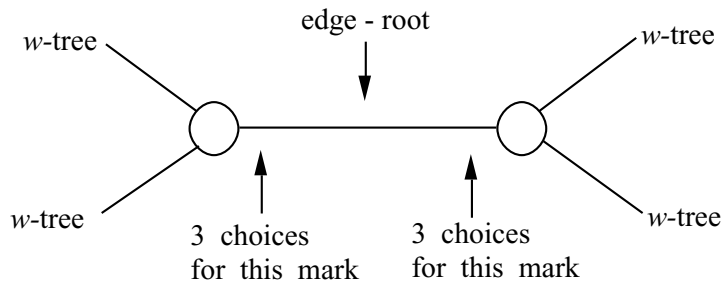}}
\caption{Edge-rooted tree} \label{Fi:edgeroot}
\end{center}
\end{figure}

Next we count  versions of these trees rooted at an
edge joining two trivalent vertices.
With the help of Figure \ref{Fi:edgeroot} we see that the
exponential generating function for such edge-rooted trees is
$(3w^2 z)^2/2 = 9w^4z^2\!/2$.

An unrooted tree with $n$ trivalent vertices is counted $n$ times
in the generating function for vertex-rooted trees and  $n-1$ times in the generating function for edge-rooted trees. Thus the difference 
$w^3z - 9w^4z^2\!/2$ counts every unrooted tree once. 
A straightforward computation shows that $w^3z - 9w^4z^2\!/2$ can be expressed most simply as $(w-u)(3u-w)/6$.
Thus the exponential generating function for graphs
whose connected components are trees is 
\begin{equation}
\label{e:nocycles}
e^{(w-u)(3u-w)/6}.
\end{equation}

As pointed out by the referee, other derivations of this generating function can be obtained by using the combinatorial interpretation of derivatives. Here is one such derivation: If $T$ is the generating function for these trees, then 
\begin{equation}
\label{e:referee}
\frac {dT}{du} = w-u = \frac{1-\sqrt{1-12uz}}{6z}-u
\end{equation} 
since if we remove one of the monovalent vertices (weighted $u$) from a tree counted by $T$, leaving an unmatched half-edge, we obtain a $w$-tree with at least one trivalent vertex. Integrating with respect to $u$,  and using the fact that $T$ has constant term 0,
we get 
\begin{equation*}
T = \frac{(1-12uz)^{3/2} -1}{108z^2}   +\frac u{6z} -\frac {u^2}2,
\end{equation*}
and this is easily checked to be equal to $(w-u)(3u-w)/6$.

We also note that the coefficients of $T$ are given explicitly by 
$$
T=\sum_{n=1}^\infty 3^n \frac{(2n)!}{(n+2)!} u^{n+2}\frac{z^n}{n!}.
$$

\begin{figure}[htbp]
\begin{center}
\scalebox{1.2}{\includegraphics{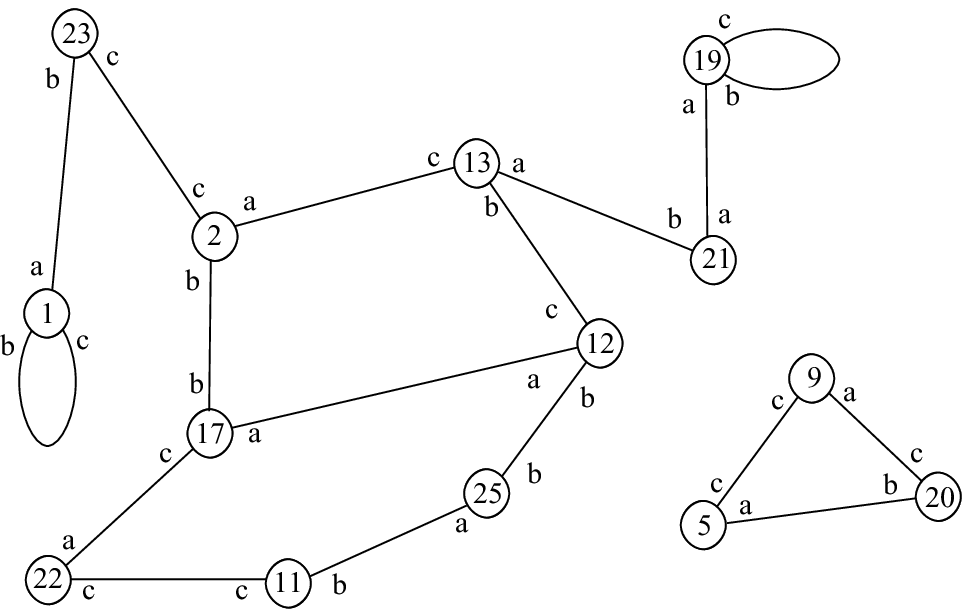}}
\caption{Components with cycles with $w$-trees eliminated}
\label{Fi:treeselim}
\end{center}
\end{figure}

\begin{figure}[htbp]
\begin{center}
\scalebox{.65}{\includegraphics{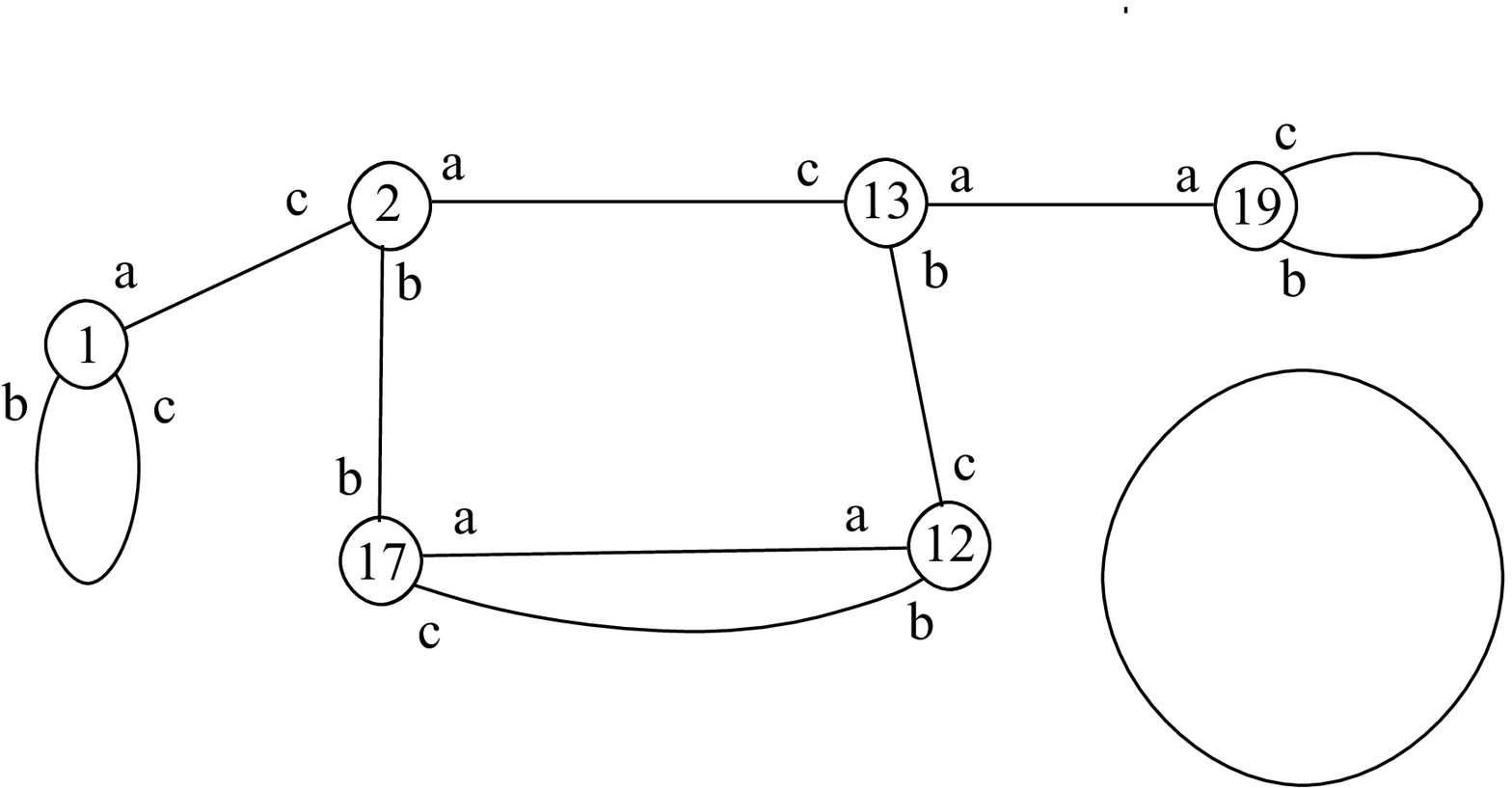}} \caption{Graph
in Figure \ref{Fi:treeselim} with bivalent vertices eliminated }
\label{Fi:finalgraph}
\end{center}
\end{figure}

Now let us look at the connected components with one or more cycles. The two
components other than the tree in  Figure
\ref{Fi:bivertelim} are of this form. First we reduce these components by
shrinking all the $w$-trees present in them. Figure
\ref{Fi:treeselim} illustrates this process.
As we see in the picture, the process of shrinking $w$-trees
leaves behind vertices weighted with $z$ which are now bivalent.
We further reduce the components by eliminating these newly
created bivalent vertices. Figure \ref{Fi:finalgraph} shows what
we get after this further reduction.
Thus we get either circles with no vertices or components (with
multiple edges and loops) in which each vertex is trivalent. We
will consider the two cases separately.

\subsection{Components with one cycle.} \label{onecyclecomp}
If the reduction of a component of a graph in $G$
results in a circle, then the original component
has exactly one cycle and each trivalent vertex on the
cycle has a $w$-tree attached at its third half-edge.
The exponential generating
function for directed cycles is
$\sum_{n=1}^{\infty}z^n/n$. To count directed cycles with $w$-trees attached,
we replace $z$ by $wz$,  and to add three marks to the half-edges of each
vertex in the cycle, we multiply 
$(wz)^n$ by~$6^n$. Finally, to undirect the cycles, we divide by 2. (Note that
marking the half-edges for the cases $n=1$ and 2 destroys  the symmetry that
would prevent us from dividing by 2 before the marks are added.)

Thus the exponential generating function for cycles
of $w$-trees is
$$\frac{1}{2}\sum_{n=1}^{\infty}
\frac{{(6wz)}^n}{n}
=\tfrac12\log{(1-6wz)}$$ and so the exponential generating function for
graphs whose connected components are cycles of $w$-trees is

\begin{equation}
\label{e:onecycle}
e^{\frac12\log {(1-6wz)}}=\frac1{\sqrt{1-6wz}}.
\end{equation}

\subsection{Components with more than one cycle.}

The components of graphs in $G$ with more than one cycle reduce to connected
graphs with only trivalent vertices. The connected graph on the left in Figure
\ref{Fi:finalgraph} is of this type. Rather than counting connected graphs of
this type, as we did in the previous two sections, and then exponentiating, we
count graphs, not necessarily connected, whose connected components are of this
type. These graphs with $m$ labeled vertices are precisely the graphs counted
by $h_{3m}(u)$ with no monovalent vertices. To see this, recall that in section~\ref{S:leftside} we described how each labeled vertex has three marked
vertices attached to it and after constructing the matching on the marked
vertices, the bivalent (or matched) marked vertices are eliminated as a
preliminary simplification. If there are no unmatched marked vertices, the
preliminary simplification gives a graph with only trivalent labeled vertices
which have marks on their half-edges.  
Thus these graphs are counted by $h_{3m}(0)$, the number of complete matchings
of $3m$ vertices. This number is 0 for $m$ odd and $(6n)!/2^{3n} (3n)!$ for
$m=2n$.  

We can recover the original graph from the reduced graph by introducing an
ordered sequence of trivalent vertices on each edge, giving each new vertex
marked half-edges and attaching a $w$-tree to one of the half-edges. We can see
this for the component with trivalent vertices in Figure~\ref{Fi:finalgraph} by
tracing it back to its original component in Figure~\ref{Fi:bivertelim}. Thus
we get a factor
$\sum_{k=0}^{\infty} {(6wz)}^k=1/(1-6wz)$ for each of the $3n$ edges in the
reduced graph, where the 6 is the number of ways to specify the marks
at a trivalent vertex. Hence the exponential generating function for graphs whose
components have more than one cycle is
\begin{equation}
\label{e:morethanonecycle}
\S \frac{(6n)!}{2^{3n}(3n)!}\frac1{(1-6wz)^{3n}} \frac{z^{2n}}{(2n)!}.
\end{equation}

By the product formula for exponential generating functions, we multiply 
\eqref{e:nocycles}, \eqref{e:onecycle}, and \eqref{e:morethanonecycle}
to get
\[
e^{(w-u)(3u-w)/6 } \frac {1} {\sqrt{1-6wz}}  \S
\frac{(6n)!}{2^{3n}(3n)!(1-6wz)^{3n}}
\frac{z^{2n}}{(2n)!},
\]
as the exponential generating function for all graphs in $G$, and  this is the
right side of~\eqref{E:mainform}.

\let\section=\oldsection 


\end{document}